# THE ASYMPTOTIC DISTRIBUTION AND BERRY–ESSEEN BOUND OF A NEW TEST FOR INDEPENDENCE IN HIGH DIMENSION WITH AN APPLICATION TO STOCHASTIC OPTIMIZATION

By Wei-Dong Liu,[1] Zhengyan Lin[1] and Qi-Man Shao[2]

*Zhejiang University, Zhejiang University and HKUST*

Let $\mathbf{X}_1,\dots,\mathbf{X}_n$ be a random sample from a $p$-dimensional population distribution. Assume that $c_1 n^\alpha \le p \le c_2 n^\alpha$ for some positive constants $c_1, c_2$ and $\alpha$. In this paper we introduce a new statistic for testing independence of the $p$-variates of the population and prove that the limiting distribution is the extreme distribution of type I with a rate of convergence $O((\log n)^{5/2}/\sqrt{n})$. This is much faster than $O(1/\log n)$, a typical convergence rate for this type of extreme distribution. A simulation study and application to stochastic optimization are discussed.

**1. Introduction and main results.** Consider a $p$-variable population represented by a random vector $\mathbf{X} = (X_1,\dots,X_p)$ with the covariance matrix $\Sigma$ and the correlation coefficient matrix R and let $\{\mathbf{X}_1,\dots,\mathbf{X}_n\}$ be a random sample of size $n$ from the population. In the applications of multivariate analysis in the century of data, both the dimension $p$ and the sample size $n$ can be very large, $p$ may be comparable with $n$ or even much larger than $n$; see, for example, Donoho [4], and Fan and Li [11]. Since classical limit theorems for fixed $p$ may not be valid for large $p$, it is necessary to develop new limiting theorems for large $p$. Whether the $p$-variables are independent is usually a primary step because independence seems to be a granted assumption for many limiting theorems. When $n/p \to \gamma > 0$ and the population distribution is normal, several statistics have been developed to test

Received September 2007; revised February 2008.

[1]Supported by National Natural Science Foundations of China (10571159, 10671176) and Specialized Research Fund for the Doctor Program of Higher Education (20060335032).

[2]Supported in part by Hong Kong RGC 602206.

the complete independence of the $p$ components of $\mathbf{X}$. Johnstone [15] uses the largest eigenvalue of the sample covariance matrix and Ledoit and Wolf [17] use quadratic forms of the sample covariance matrix eigenvalues to test the null hypothesis $H_0 : \Sigma = I_p$, where $I_p$ is the $p \times p$ identity matrix, while Schott [20] uses the sums of sample correlation coefficient squares to test $H_0 : \mathrm{R} = I_p$. When the normality is not present, Jiang [14] constructs a test statistic based on the largest entries of the sample correlation matrix. Write

$$\mathbf{X}_k = (X_{k,1}, X_{k,2}, \ldots, X_{k,p}), \qquad 1 \le k \le n$$

and let

$$\widetilde{L}_n = \max_{1 \le i < j \le p} |\widetilde{\rho}_{i,j}^{(n)}|,$$

where

$$\widetilde{\rho}_{i,j}^{(n)} = \frac{\sum_{k=1}^{n}(X_{k,i} - \bar{X}_i^{(n)})(X_{k,j} - \bar{X}_j^{(n)})}{(\sum_{k=1}^{n}(X_{k,i} - \bar{X}_i^{(n)})^2)^{1/2}(\sum_{k=1}^{n}(X_{k,j} - \bar{X}_j^{(n)})^2)^{1/2}}$$

and $\bar{X}_i^{(n)} = \sum_{k=1}^{n} X_{k,i}/n$. Jiang [14] proves the following limit theorem concerning the test statistic $\widetilde{L}_n$:

*If $n/p \to \gamma \in (0, \infty)$ and $\mathsf{E}|X_{11}|^r < \infty$ for some $r > 30$, then*

$$\lim_{n \to \infty} \mathsf{P}(n\widetilde{L}_n^2 - 4\log p + \log_2 p \le y) = \exp(-e^{-y/2}/\sqrt{8\pi}) \tag{1.1}$$

*for $y \in R$, where and in the sequel $\log x = \ln \max(x, e)$ and $\log_2 x = \log(\log x)$.*

Zhou [21] shows that the moment condition $\mathsf{E}|X_{1,1}|^r < \infty$ for some $r > 30$ can be weakened to

$$x^6 \mathsf{P}(|X_{1,1}X_{1,2}| \ge x) \to 0 \qquad \text{as } x \to \infty. \tag{1.2}$$

The limit distribution appearing in (1.1) is called the extreme distribution of type I. It seems a common belief that the convergence rate of this type of extreme distribution is typically slow (see Hall [13]). In fact, we shall prove [see Theorem 1.2 and (1.11)] that even when $X_{1,1}$ has the standard normal distribution the rate of convergence is of order of $O(\log_2 n/\log n)$. The main purpose of this paper is to introduce a modified test statistic and show that the new one also has an extreme limit distribution of type I, but with a rate of convergence of $O((\log n)^{5/2}/\sqrt{n})$. We shall also prove that the approximation rate of $\mathsf{P}(n\widetilde{L}_n^2 - 4\log p + \log_2 p \le y)$ to $\exp(-\frac{p^2-p}{2}\mathsf{P}(\chi^2(1) \ge 4\log p - \log_2 p + y))$ instead of the final limit $\exp(-e^{-y/2}/\sqrt{8\pi})$ is indeed of order of $O((\log n)^{5/2}/\sqrt{n})$. This indicates that when a statistic has an extreme limiting distribution, one should use some "intermediate" approximation, not the final limiting distribution to approximate the distribution of the statistic. Extreme limiting distributions are important in various applications, including assessing risk for highly unusual events, hydrologic assessment, analysis of network simulation and engineering (see [12] and [16]). These results are partially motivated by the need for new approaches to applications of these types in practice.

Throughout this paper let $H_0$ be the null hypothesis that the $p$ components of the population $\mathbf{X}$ are independent and have the same distribution, and let $\mathbf{X}_k = (X_{k,1}, X_{k,2}, \dots, X_{k,p}), 1 \le k \le n$, be a random sample from the population $\mathbf{X}$. Define

$$L_n^2 = \max_{1 \le i < j \le p} r_{i,j}^2, \tag{1.3}$$

where

$$\begin{aligned} r_{i,j}^2 &= (2A_{n,i,j}^2 + 2B_{n,i,j}^2)/D_{n,i,j}, \\ A_{n,i,j} &= \sum_{k=1}^{[n/2]} (X_{k,i} - \bar{X}_i^{(n)})(X_{k,j} - \bar{X}_j^{(n)}), \\ B_{n,i,j} &= \sum_{k=[n/2]+1}^{n} (X_{k,i} - \bar{X}_i^{(n)})(X_{k,j} - \bar{X}_j^{(n)}), \\ D_{n,i,j} &= \sum_{k=1}^{n} (X_{k,i} - \bar{X}_i^{(n)})^2 \sum_{k=1}^{n} (X_{k,j} - \bar{X}_j^{(n)})^2, \end{aligned} \tag{1.4}$$

where $[n/2]$ denotes the integer part of $n/2$. Our test statistic is chosen as

$$\mathcal{W}_n := nL_n^2 - 4\log p. \tag{1.5}$$

It is easy to see that $r_{i,j}^2$ is a consistent estimator of the square of the correlation coefficient between $X_i$ and $X_j$ and that $r_{i,j}^2 \ge (\tilde{\rho}_{i,j}^{(n)})^2$.

Instead of assuming that $n$ and $p$ have the same order, we consider a more general case. Assume

$$c_1 n^\alpha \le p \le c_2 n^\alpha, \tag{1.6}$$

where $c_1, c_2$ and $\alpha$ are positive constants.

Our first result shows that the $\mathcal{W}_n$ has an extreme limiting distribution of type I under a weaker moment assumption than (1.2) and the error of the approximation is of order of $O((\log n)^{5/2}/\sqrt{n})$ under $\mathsf{E}|X_{1,1}|^{3+4\alpha} < \infty$ and $\alpha > 3/4$.

THEOREM 1.1. *Assume that (1.6) is satisfied and that*

$$x^{1+2\alpha} \mathsf{P}(|X_{1,1}X_{1,2}| \ge \sqrt{x \log x}) \to 0 \qquad \textit{as } x \to \infty. \tag{1.7}$$

*Then under* $H_0$,

$$\mathsf{P}(\mathcal{W}_n \le y) \to \exp\biggl(-\frac{1}{2}\exp\biggl(-\frac{y}{2}\biggr)\biggr) \tag{1.8}$$

*as* $n \to \infty$ *for any* $y \in R$. *If* $\mathsf{E}|X_{1,1}|^{3+4\alpha} < \infty$ *and* $\alpha > 3/4$, *then*

$$\sup_{y\in R}\biggl|\mathsf{P}(\mathcal{W}_n \le y) - \exp\biggl(-\frac{1}{2}\exp\biggl(-\frac{y}{2}\biggr)\biggr)\biggr| \le Cn^{-1/2}(\log n)^{5/2} \tag{1.9}$$

*and* $C$ *is a constant independent of* $n$ *and* $p$.

The next result shows that (1.1) remains valid under the assumption (1.7). Moreover, the convergence rate to $\exp(-\frac{p^2-p}{2}\mathsf{P}(\chi^2(1) \ge 4\log p - \log_2 p + y))$ [$\chi^2(1)$ has a chi-square distribution with 1 degree of freedom] can also achieve the order of $O((\log n)^{5/2}/\sqrt{n})$ when the 7th moment of $X_{1,1}$ is finite.

THEOREM 1.2. *Under* $H_0$, *if (1.6) and (1.7) are satisfied, then (1.1) holds. If* $\mathsf{E}|X_{1,1}|^{3+4\alpha} < \infty$ *and* $\alpha > 3/4$, *then*

$$\begin{aligned}\sup_{y\in R}\biggl|&\mathsf{P}(n\widetilde{L}_n^2 - 4\log p + \log_2 p \le y)\\ &- \exp\biggl(-\frac{p^2-p}{2}\mathsf{P}(\chi^2(1) \ge 4\log p - \log_2 p + y)\biggr)\biggr|\\ &\le Cn^{-1/2}(\log n)^{5/2}.\end{aligned} \tag{1.10}$$

One can show that (see a proof in Section 5) if $\alpha = 1$ in (1.6), then

$$\begin{aligned}&\exp\biggl(-\frac{p^2-p}{2}\mathsf{P}(\chi^2(1) \ge 4\log p - \log_2 p + y)\biggr) - \exp(-e^{-y/2}/\sqrt{8\pi})\\ &\quad\sim \frac{\log_2 n}{8\log n}\frac{1}{\sqrt{8\pi}}\exp\biggl(-\frac{y}{2} - \frac{1}{\sqrt{8\pi}}\exp\biggl(-\frac{y}{2}\biggr)\biggr).\end{aligned} \tag{1.11}$$

So the rate of convergence in (1.1) is of order of $O(\log_2 n/\log n)$.

The following remarks are noted.

REMARK 1.1. The logarithmic term $(\log n)^{5/2}$ in (1.10) may not be sharp. Since our aim is to get the main order $n^{-1/2}$, we will not seek the optimality of the logarithmic term in this paper.

REMARK 1.2. It is not necessary to require $p$ and $n$ have the tight relation (1.6). For example, when (1.6) is replaced by

$$c_1 n^{\alpha_1} \le p \le c_2 n^{\alpha}$$

for some positive constants $c_1, c_2, \alpha_1$ and $\alpha$, and assume that $\mathsf{E}|X_{1,1}|^{(3+4\alpha)\vee r} < \infty$ for some $r > 6$, then, following the proofs of Theorems 1.1 and 1.2, we have

$$\sup_{y\in R}\left|\mathsf{P}(\mathcal{W}_n \le y) - \exp\left(-\frac{1}{2}\exp\left(-\frac{y}{2}\right)\right)\right| \le C(n^{-1/2}(\log p)^{5/2} + p^{-1+\varepsilon})$$

and

$$\begin{aligned}\sup_{y\in R}\Big|&\mathsf{P}(n\widetilde{L}_n^2 - 4\log p + \log_2 p \le y)\\ &- \exp\left(-\frac{p^2-p}{2}\mathsf{P}(\chi^2(1) \ge 4\log p - \log_2 p + y)\right)\Big|\\ &\le C(n^{-1/2}(\log p)^{5/2} + p^{-1+\varepsilon})\end{aligned}$$

for $\varepsilon > 0$, where the constant $C$ may depend on $\varepsilon$. When $\log p$ is of order of $n^r$ $(0 < r \le 1)$, we shall discuss the case in our future paper.

REMARK 1.3. It is not necessary to assume that the $p$ components of $\mathbf{X}$ are identically distributed in Theorems 1.1 and 1.2. For example, if $\inf_n \min_{1\le j\le p} \mathsf{E}X_{1,j}^2 > 0$ and

$$x^{1+2\alpha}\sup_{n\ge 1}\max_{1\le i<j\le p}\mathsf{P}(|X_{1,i}X_{1,j}| \ge \sqrt{x\log x}) \to 0,$$

then (1.8) remains valid. If, in addition, $\sup_n \max_{1\le j\le p} \mathsf{E}|X_{1,j}|^{3+4\alpha} < \infty$, $\alpha > 3/4$, then (1.9) and (1.10) hold.

The paper is organized as follows. In Section 2 we conduct a simulation study and give a simple application of Theorem 1.2 to the sparsest solution of large underdetermined systems of linear equations. In Section 3 we present a general theorem from which Theorems 1.1 and 1.2 can be derived easily. An outline of the proof of the general result along with five propositions is given in Section 5, while the detailed proofs of the propositions are postponed to Section 6. Proofs of Theorems 1.1 and 1.2 are given in Section 4.

Throughout the paper, $C$ will denote a positive constant that doesn't depend on $n$ and $p$ but may be different at each appearance.

**2. A simulation study and application to stochastic optimization.** In this section we give a simulation study for performance of $\mathcal{W}_n$ and $\tilde{L}_n$ and an application to the sparsest solution of large underdetermined systems of linear equations.

TABLE 1

*Estimated significance levels when $\alpha = 0.05$ and $X_{1,1} \sim N(0,1)$*

| $p$ | Test statistics | $n = 16$ | $n = 32$ | $n = 64$ | $n = 128$ | $n = 256$ |
|---|---|---|---|---|---|---|
| 4 | $\mathcal{W}_n$ | 0.0484 | 0.0580 | 0.0496 | 0.0520 | 0.0566 |
| | $\tilde{L}_{\text{new}}$ | 0.0318 | 0.0458 | 0.0498 | 0.0474 | 0.0532 |
| | $\tilde{L}_{\text{old}}$ | 0.0140 | 0.0232 | 0.0284 | 0.0256 | 0.0292 |
| 8 | $\mathcal{W}_n$ | 0.0190 | 0.0332 | 0.0412 | 0.0524 | 0.0440 |
| | $\tilde{L}_{\text{new}}$ | 0.0104 | 0.0316 | 0.0368 | 0.0462 | 0.0462 |
| | $\tilde{L}_{\text{old}}$ | 0.0066 | 0.0198 | 0.0222 | 0.0258 | 0.0338 |
| 16 | $\mathcal{W}_n$ | 0.0094 | 0.0248 | 0.0356 | 0.0436 | 0.0478 |
| | $\tilde{L}_{\text{new}}$ | 0.0012 | 0.0112 | 0.0316 | 0.0420 | 0.0482 |
| | $\tilde{L}_{\text{old}}$ | 0.0002 | 0.0130 | 0.0246 | 0.0312 | 0.0338 |
| 32 | $\mathcal{W}_n$ | 0.0032 | 0.0188 | 0.0366 | 0.0412 | 0.0432 |
| | $\tilde{L}_{\text{new}}$ | 0.0000 | 0.0094 | 0.0228 | 0.0368 | 0.0376 |
| | $\tilde{L}_{\text{old}}$ | 0.0000 | 0.0044 | 0.0212 | 0.0280 | 0.0364 |
| 64 | $\mathcal{W}_n$ | 0.0010 | 0.0114 | 0.0296 | 0.0402 | 0.0460 |
| | $\tilde{L}_{\text{new}}$ | 0.0000 | 0.0020 | 0.0100 | 0.0292 | 0.0356 |
| | $\tilde{L}_{\text{old}}$ | 0.0000 | 0.0026 | 0.0160 | 0.0256 | 0.0358 |
| 128 | $\mathcal{W}_n$ | 0.0004 | 0.0082 | 0.0218 | 0.0350 | 0.0568 |
| | $\tilde{L}_{\text{new}}$ | 0.0000 | 0.0000 | 0.0060 | 0.0262 | 0.0380 |
| | $\tilde{L}_{\text{old}}$ | 0.0000 | 0.0000 | 0.0060 | 0.0170 | 0.0362 |

2.1. *A simulation study.* The performance of the test statistics $\mathcal{W}_n$ and $\tilde{L}_n$ is carried out via simulation. Estimates of the actual significance levels are obtained from 5000 independent simulations with the nominal significance level $\alpha = 0.05$.

The simulation results for tests of $H_0$ based on $\mathcal{W}_n$ in Theorem 1.1, $\tilde{L}_n$ (denoted by $\tilde{L}_{\text{new}}$) in Theorem 1.2 and $\tilde{L}_n$ (denoted by $\tilde{L}_{\text{old}}$) in (1.1) are given in Table 1 when $X_{1,1}$ has a standard normal distribution, and in Table 2 when $X_{1,1}$ has a $t$-distribution with 7 degrees of freedom. The estimated significance levels are usually lower than the nominal level 0.05, which indicates that the tests are conservative. The performances of $W_n$ and $\tilde{L}_{\text{new}}$ are comparable and both are well when $n$ is larger than $p$.

2.2. *An application to stochastic optimization.* We apply Theorem 1.2 to the problem of finding sparse representations of single measurement vectors (SMV) in an over-complete dictionary. The SMV problem can be described as follows. Given a single measurement vector $b$ and a dictionary $A$, one wants to solve the system of equations $Ax = b$, where $A$ is a $n \times p$ matrix, $x$ is a $p$-variable vector and $b$ is a $n$-variable vector. It is usually assumed that $n \ll p$. A sparse representation means that vector $x$ has a small number of nonzero components. Examples of such underdetermined systems of equations include array signal processing, inverse problems and genomic data analysis. We refer to [2, 6, 7, 8] and [9] and references therein for a comprehensive description of many important applications of a SMV problem.

TABLE 2
*Estimated significance levels when $\alpha = 0.05$ and $X_{1,1} \sim t_7$*

| $p$ | Test statistics | $n = 16$ | $n = 32$ | $n = 64$ | $n = 128$ | $n = 256$ |
|---|---|---|---|---|---|---|
| 4 | $\mathcal{W}_n$ | 0.0484 | 0.0580 | 0.0496 | 0.0520 | 0.0566 |
| | $\tilde{L}_{\text{new}}$ | 0.0408 | 0.0466 | 0.0514 | 0.0504 | 0.0522 |
| | $\tilde{L}_{\text{old}}$ | 0.0192 | 0.0242 | 0.0286 | 0.0276 | 0.0300 |
| 8 | $\mathcal{W}_n$ | 0.0374 | 0.0560 | 0.0624 | 0.0574 | 0.0586 |
| | $\tilde{L}_{\text{new}}$ | 0.0084 | 0.0332 | 0.0468 | 0.0440 | 0.0470 |
| | $\tilde{L}_{\text{old}}$ | 0.0054 | 0.0230 | 0.0290 | 0.0358 | 0.0372 |
| 16 | $\mathcal{W}_n$ | 0.0292 | 0.0536 | 0.0750 | 0.0676 | 0.0622 |
| | $\tilde{L}_{\text{new}}$ | 0.0012 | 0.0186 | 0.0324 | 0.0436 | 0.0446 |
| | $\tilde{L}_{\text{old}}$ | 0.0004 | 0.0146 | 0.0308 | 0.0336 | 0.0366 |
| 32 | $\mathcal{W}_n$ | 0.0144 | 0.0682 | 0.0886 | 0.0758 | 0.0664 |
| | $\tilde{L}_{\text{new}}$ | 0.0000 | 0.0102 | 0.0294 | 0.0414 | 0.0444 |
| | $\tilde{L}_{\text{old}}$ | 0.0000 | 0.0062 | 0.0244 | 0.0336 | 0.0436 |
| 64 | $\mathcal{W}_n$ | 0.0066 | 0.0816 | 0.1122 | 0.1010 | 0.0670 |
| | $\tilde{L}_{\text{new}}$ | 0.0000 | 0.0040 | 0.0266 | 0.0382 | 0.0472 |
| | $\tilde{L}_{\text{old}}$ | 0.0000 | 0.0042 | 0.0196 | 0.0352 | 0.0388 |
| 128 | $\mathcal{W}_n$ | 0.0000 | 0.1010 | 0.1240 | 0.1138 | 0.0820 |
| | $\tilde{L}_{\text{new}}$ | 0.0000 | 0.0002 | 0.0184 | 0.0354 | 0.0480 |
| | $\tilde{L}_{\text{old}}$ | 0.0000 | 0.0002 | 0.0178 | 0.0342 | 0.0438 |

A sparse representation can be found by solving the following optimization problem:

$$\text{(Q0)} \qquad \min \|x\|_0, \qquad \text{s.t. } Ax = b,$$

where the quantity $\|x\|_0$ denotes the number of nonzero elements in the vector $x$. The problem (Q0) is essentially a combinatorial optimization problem, which, in general, is extremely difficult to solve. The above problem can be convexified as a $\ell_1$-norm minimization problem, and solved via linear programming. The $\ell_1$-norm minimization problem is

$$\text{(Q1)} \qquad \min \|x\|_1, \qquad \text{s.t. } Ax = b,$$

where $\|x\|_1$ is the sum of the absolute values of the elements of vector $x$. It has been proved that the solutions between (Q0) and (Q1) are equivalent under various conditions. For example, letting $G = A^T A$ and $M = \max_{1 \le i,j \le p, i \ne j} |G(i,j)|$, if $\|x\|_0 < (1 + M^{-1})/2$, then $x$ is the unique solution of (Q1) (for $b = Ax$) and this solution is identical to the unique solution

of (Q0) [3]. Here we assume that $A$ is a random matrix. Let $n, p$ satisfy the condition in Remark 1.2 and $\{X_{k,i}, k, i \geq 1\}$ be independent centered random variables satisfying the conditions in Remark 1.3. Define the normalized $(k,i)$ element of $A$ by $Y_{k,i} := X_{k,i}/(\sum_{k=1}^n X_{k,i}^2)^{1/2}$. Following the proof of Theorem 1.2, we see that (1.10) remains valid for $M^2 = M_n^2$, where

$$M_n^2 = \max_{1\leq i,j\leq p, i\neq j} |G_{i,j}|^2 = \max_{1\leq i<j\leq p} \frac{(\sum_{k=1}^n X_{k,i}X_{k,j})^2}{(\sum_{k=1}^n X_{k,i}^2)(\sum_{k=1}^n X_{k,j}^2)}.$$

Hence, (Q0) and (Q1) are equivalent with probability $1-\alpha$ $(0<\alpha<1)$ for every $x$ with fewer than $(1+m_\alpha)/2$ nonzeros, where $m_\alpha = \sqrt{n/(y_\alpha + 4\log p - \log_2 p)}$ and $y_\alpha$ is the solution of

$$\exp\biggl(-\frac{p^2-p}{2}\mathsf{P}(\chi^2(1) \geq 4\log p - \log_2 p + y_\alpha)\biggr) = 1-\alpha.$$

When $X_{ij}$ above are i.i.d. standard normal random variables, the result is similar to that given in [5].

**3. A general result.** Instead of proving Theorems 1.1 and 1.2 separately, we give a general result in this section. Let $d$ be a positive integer and $\{X, X_{k,i}^{(m)}; k, i \geq 1, 1 \leq m \leq d\}$ be an array of i.i.d. random variables. Put

$$X_{k,i,j} = (Y_{k,i,j}^{(1)}, \ldots, Y_{k,i,j}^{(d)}), \qquad Y_{k,i,j}^{(m)} = X_{k,i}^{(m)} X_{k,j}^{(m)}, \qquad i,j,k \geq 1, 1 \leq m \leq d$$

and

$$W_{p,n} = \max_{1\leq i<j\leq p} \Biggl\|\sum_{k=1}^n X_{k,i,j}\Biggr\|,$$

where $\|\cdot\|$ denotes the Euclidean norm in $R^d$.

THEOREM 3.1. *Suppose that $\mathsf{E}X = 0$ and $\mathsf{E}X^2 = 1$. Let $X'$ be an independent copy of $X$, and*

$$\sup_x x^{1+2\alpha}\mathsf{P}(|XX'| \geq \sqrt{x\log x}) < \infty. \tag{3.1}$$

*Then for any $0 < \varepsilon \leq 10^{-4}$ there exists a finite constant $C$ such that*

$$\begin{aligned}\sup_{y\in R}&\biggl|\mathsf{P}\biggl(\frac{W_{p,n}^2}{n} - \alpha_p \leq y\biggr) - \exp\biggl(-\frac{p^2-p}{2}\mathsf{P}(\chi^2(d) \geq \alpha_p + y)\biggr)\biggr| \\ &\leq Cp^{-1+20\sqrt{\varepsilon}} + C\frac{(\log n)^{5/2}}{n^{1/2}}\mathsf{E}|XX'|^3 I\{|XX'| \leq \sqrt{n}/(\log n)^4\} \\ &\quad + Cn^{1+2\alpha}\mathsf{P}(|XX'| \geq d^{-1/2}\varepsilon\sqrt{n\log p}),\end{aligned} \tag{3.2}$$

*where $\alpha_p = 4\log p - (2-d)\log_2 p$ and $\chi^2(d)$ has a chi-square distribution with $d$ degrees of freedom.*

Now, set

$$A_{n,i} = \sum_{m=1}^{d} \sum_{k=1}^{n} (X_{k,i}^{(m)})^2, \qquad 1 \le i \le p$$

and $Q_{n,i,j} = A_{n,i} A_{n,j}$,

$$\mathcal{L}_{p,n}^2 = \max_{1 \le i < j \le p} \frac{\| \sum_{k=1}^{n} X_{k,i,j} \|^2}{Q_{n,i,j}}.$$

THEOREM 3.2. *Under the conditions of Theorem 3.1 and* $\mathsf{E}X^4 < \infty$, *we have, for any* $0 < \varepsilon \le 10^{-4}$,

$$\sup_{y \in R} \left| \mathsf{P}(d^2 n \mathcal{L}_{p,n}^2 - \alpha_p \le y) - \exp\left( -\frac{p^2 - p}{2} \mathsf{P}(\chi^2(d) \ge \alpha_p + y) \right) \right|$$

$$\begin{aligned} &\le C p^{-1+20\sqrt{\varepsilon}} + C \frac{(\log n)^{5/2}}{n^{1/2}} \\ &\quad + C n^{1+2\alpha} \mathsf{P}(|XX'| \ge d^{-1/2} \varepsilon \sqrt{n \log p}) + \tau_n, \end{aligned} \tag{3.3}$$

*where*

$$\tau_n = C n p^{20\sqrt{\varepsilon}} \mathsf{P}\left( |X| \ge \frac{n^{1/4}}{(\log p)^{1/4}} \right).$$

The proofs of Theorems 3.1 and 3.2 are postponed to Section 5.

**4. Proofs of Theorems 1.1 and 1.2.** We are now ready to see that Theorems 1.1 and 1.2 are two special cases of Theorems 3.1 and 3.2.

PROOF OF THEOREM 1.1. For the sake of simplicity, we assume that $n$ is even. Otherwise use $[n/2]$ instead of $n/2$ below. Without loss of generality, we assume $\mathsf{E}X_{1,1} = 0$ and $\mathsf{E}X_{1,1}^2 = 1$. Set

$$\widetilde{A}_{n,i,j} = \sum_{k=1}^{n/2} X_{k,i} X_{k,j}, \qquad \widetilde{B}_{n,i,j} = \sum_{k=n/2+1}^{n} X_{k,i} X_{k,j} = \sum_{k=1}^{n/2} X_{k+n/2,i} X_{k+n/2,j},$$

$$\widetilde{D}_{n,i,j} = \sum_{k=1}^{n} X_{k,i}^2 \sum_{k=1}^{n} X_{k,j}^2.$$

Take $d = 2$ in Theorem 3.1. Since $(\widetilde{A}_{n,i,j})^2 + (\widetilde{B}_{n,i,j})^2 = \|(\widetilde{A}_{n,i,j}, \widetilde{B}_{n,i,j})\|^2$, by Theorem 3.1,

$$\begin{aligned} &\mathsf{P}\left( \max_{1 \le i < j \le p} \frac{2(\widetilde{A}_{n,i,j})^2 + 2(\widetilde{B}_{n,i,j})^2}{n} - 4 \log p \le y \right) \\ &\quad = \mathsf{P}\left( \frac{2W_{p,n/2}^2}{n} - 4 \log p \le y \right) \to e^{-e^{-y/2}/2} \end{aligned} \tag{4.1}$$

and hence,

$$\max_{1\le i<j\le p}|\widetilde{A}_{n,i,j}| = O_{\mathsf{P}}(\sqrt{n\log n}) \quad \text{and} \tag{4.2}$$
$$\max_{1\le i<j\le p}|\widetilde{B}_{n,i,j}| = O_{\mathsf{P}}(\sqrt{n\log n}).$$

Noting that condition (1.7) implies $\mathsf{E}|X_{1,1}|^{2+4\alpha-\varepsilon}<\infty$ for any $\varepsilon>0$ (see Lemma 6.4), we can show, by Lemma 6.1,

$$\max_{1\le i\le p}\frac{|\sum_{k=1}^n X_{k,i}|}{n} = O_{\mathsf{P}}\left(\sqrt{\frac{\log n}{n}}\right) \tag{4.3}$$

and

$$\max_{1\le i\le p}\frac{|\sum_{k=1}^n X_{k,i}^2-n|}{n} = O_{\mathsf{P}}(n^{-\delta}), \tag{4.4}$$

for some $\delta>0$. Observe that

$$\begin{aligned} r_{i,j}^2 &= 2\frac{(\sum_{k=1}^{n/2}X_{k,i}X_{k,j})^2-2\bar{X}_{1,i,j}\sum_{k=1}^{n/2}X_{k,i}X_{k,j}+(\bar{X}_{1,i,j})^2}{[\sum_{k=1}^n X_{k,i}^2-n(\bar{X}_i^{(n)})^2][\sum_{k=1}^n X_{k,j}^2-n(\bar{X}_j^{(n)})^2]} \\ &\quad + 2\frac{(\sum_{k=n/2+1}^{n}X_{k,i}X_{k,j})^2-2\bar{X}_{2,i,j}\sum_{k=n/2+1}^{n}X_{k,i}X_{k,j}+(\bar{X}_{2,i,j})^2}{[\sum_{k=1}^n X_{k,i}^2-n(\bar{X}_i^{(n)})^2][\sum_{k=1}^n X_{k,j}^2-n(\bar{X}_j^{(n)})^2]}, \end{aligned}$$

where

$$\begin{aligned} \bar{X}_{1,i,j} &= 2^{-1}n(\bar{X}_i^{(n/2)}\bar{X}_j^n+\bar{X}_i^{(n)}\bar{X}_j^{n/2}-\bar{X}_i^{(n)}\bar{X}_j^n) = O_{\mathsf{P}}(\log n), \\ \bar{X}_{2,i,j} &= 2^{-1}n\left(2n^{-1}\sum_{k=n/2+1}^n X_{k,i}\bar{X}_j^n+2n^{-1}\bar{X}_i^{(n)}\sum_{k=n/2+1}^n X_{k,j}-\bar{X}_i^{(n)}\bar{X}_j^n\right) \\ &= O_{\mathsf{P}}(\log n). \end{aligned}$$

(1.8) now follows from (4.1)–(4.4).

Now we prove (1.9). Let $E_{n,i,j}=(A_{n,i,j})^2+(B_{n,i,j})^2$. We first show that $D_{n,i,j}$ in the denominator of $r_{i,j}^2$ can be replaced by $\widetilde{D}_{n,i,j}$. Observe that

$$\frac{\sum_{k=1}^n[(X_{k,i}-\bar{X}_i^{(n)})^2-X_{k,i}^2]}{n} = -\frac{(\sum_{k=1}^n X_{k,i})^2}{n^2}.$$

We have, by Lemma 6.1,

$$\begin{aligned} &\mathsf{P}\left(\max_{1\le i\le p}\left|\frac{\sum_{k=1}^n[(X_{k,i}-\bar{X}_i^{(n)})^2-X_{k,i}^2]}{n}\right|\ge 4\sqrt{\log n/n}\right) \\ &\qquad \le p\mathsf{P}\left(\left|\sum_{k=1}^n X_{k,1}\right|\ge 2n\left(\frac{\log n}{n}\right)^{1/4}\right) \end{aligned} \tag{4.5}$$

$$\leq Cn^{-1/2},$$

$$\mathsf{P}\left(\max_{1\leq i\leq p}\left|\frac{\sum_{k=1}^{n} X_{k,i}^2 - n}{n}\right| \geq \frac{1}{2}\right) \tag{4.6}$$
$$\leq Cn^{-1/2-\alpha} \leq Cn^{-1/2}$$

and

$$\mathsf{P}\left(\max_{1\leq i\leq p}\left|\sum_{k=1}^{n} X_{k,i}\right| \geq 4\sqrt{n\log n}\right) \leq Cn^{-1/2}. \tag{4.7}$$

Therefore,

$$\mathsf{P}\left(2n\max_{1\leq i<j\leq p}\frac{E_{n,i,j}}{\widetilde{D}_{n,i,j}} \leq (1-8\sqrt{\log n/n})(y+\alpha_p)\right) - Cn^{-1/2}$$
$$\leq \mathsf{P}(\mathcal{W}_n \leq y)$$
$$\leq \mathsf{P}\left(2n\max_{1\leq i<j\leq p}\frac{E_{n,i,j}}{\widetilde{D}_{n,i,j}} \leq (1+8\sqrt{\log n/n})(y+\alpha_p)\right) + Cn^{-1/2}.$$

Now write $F_{n,i,j} = (\widetilde{A}_{n,i,j}, \widetilde{B}_{n,i,j})$. Note that

$$|\|F_{n,i,j}\| - \sqrt{A_{n,i,j}^2 + B_{n,i,j}^2}| \leq \|(\bar{X}_{1,i,j}, \bar{X}_{2,i,j})\|.$$

This together with (4.7) leads to

$$\mathsf{P}\left(\sqrt{2n}\max_{1\leq i<j\leq p}\frac{\|F_{n,i,j}\|}{\widetilde{D}_{n,i,j}^{1/2}} \leq (1-8\sqrt{\log n/n})^{1/2}(y+\alpha_p)^{1/2}\right.$$
$$\left. - C\log n/\sqrt{n}\right) - Cn^{-1/2}$$
$$\leq \mathsf{P}(\mathcal{W}_n \leq y)$$
$$\leq \mathsf{P}\left(\sqrt{2n}\max_{1\leq i<j\leq p}\frac{\|F_{n,i,j}\|}{\widetilde{D}_{n,i,j}^{1/2}} \leq (1+8\sqrt{\log n/n})^{1/2}(y+\alpha_p)^{1/2}\right.$$
$$\left. + C\log n/\sqrt{n}\right) + Cn^{-1/2}.$$

Take $d=2$ in Theorem 3.2. It is easily seen that $\mathsf{P}(2n\max_{1\leq i<j\leq p}\frac{\|F_{n,i,j}\|^2}{\widetilde{D}_{n,i,j}} \leq x) = \mathsf{P}(2n\mathcal{L}_{p,n/2}^2 \leq x)$ for any $x \in R$. Write

$$l_{n\pm}(y) = \left[\left(1 \pm 8\sqrt{\frac{\log n}{n}}\right)^{1/2}(y+\alpha_p)^{1/2} \pm C\frac{\log n}{n^{1/2}}\right]^2.$$

We have

$$\begin{aligned} &\mathsf{P}(2n\mathcal{L}^2_{p,n/2} \le l_{n-}(y)) - Cn^{-1/2} \\ &\quad \le \mathsf{P}(\mathcal{W}_n \le y) \\ &\quad \le \mathsf{P}(2n\mathcal{L}^2_{p,n/2} \le l_{n+}(y)) + Cn^{-1/2}. \end{aligned} \tag{4.8}$$

By Theorem 3.2, we obtain

$$\begin{aligned} &\sup_{y\in R}\left|\mathsf{P}(2n\mathcal{L}^2_{p,n/2} \le l_{n+}(y)) - \exp\left(-\frac{1}{2}\exp\left(-\frac{y}{2}\right)\right)\right| \\ &\quad \le \sup_{y\in R}\left|\exp\left(-\frac{p^2-p}{2}\mathsf{P}(\chi^2(2) \ge l_{n+}(y))\right) - \exp\left(-\frac{1}{2}\exp\left(-\frac{y}{2}\right)\right)\right| \\ &\qquad + Cn^{-1/2}(\log n)^{5/2} + Cp^{-1+20\sqrt{\varepsilon}} \\ &\qquad + dn^{1+2\alpha}\mathsf{P}(|XX'| \ge 2^{-1}\varepsilon\sqrt{n\log p}) + \tau_n \\ &\quad =: \sup_{y\in R}\mathsf{P}_{n+}(y) + Cn^{-1/2}(\log n)^{5/2}, \end{aligned}$$

where

$$\mathsf{P}_{n+}(y) = \left|\exp\left(-\frac{p^2-p}{2}e^{-l_{n+}(y)/2}\right) - \exp\left(-\frac{1}{2}\exp\left(-\frac{y}{2}\right)\right)\right|.$$

Note that

$$\sup_{y\le -2\log_2 n^8} l_{n+}(y) - \alpha_p \le -2\log_2 n^8 + C\frac{(\log n)^{3/2}}{n^{1/2}}. \tag{4.9}$$

This implies

$$\sup_{y\le -2\log_2 n^8} \mathsf{P}_{n+}(y) \le Cn^{-3}.$$

Also, we can get

$$\sup_{y\ge 2\log n} l_{n+}(y) - \alpha_p \ge 2\log n. \tag{4.10}$$

Therefore, following the inequality $1-\exp(-\frac{1}{2}\exp(-\frac{y}{2})) \le C\exp(-\frac{y}{2})$ for $y\ge 1$,

$$\sup_{y\ge 2\log n} \mathsf{P}_{n+}(y) \le Cn^{-1/2}.$$

A direct elementary calculation gives

$$\sup_{-2\log_2 n^8\le y\le 2\log n} |l_{n+}(y) - \alpha_p - y| \le C\frac{(\log n)^{3/2}}{n^{1/2}}, \tag{4.11}$$

hence, by the inequality $|e^x - 1| \le C|x|$ for $|x| \le 1$,

$$\sup_{-2\log_2 n^8 \le y \le 2\log n} \mathsf{P}_{n+}(y) \le Cn^{-1/2}(\log n)^{5/2} + Cp^{-1}\log n.$$

The above arguments yield

$$\sup_{y\in R}\Big|\mathsf{P}(2n\mathcal{L}^2_{p,n/2} \le l_{n+}(y)) - \exp\Big(-\frac{1}{2}\exp\Big(-\frac{y}{2}\Big)\Big)\Big| \le Cn^{-1/2}(\log n)^{5/2}.$$

Similarly,

$$\sup_{y\in R}\Big|\mathsf{P}(2n\mathcal{L}^2_{p,n/2} \le l_{n-}(y)) - \exp\Big(-\frac{1}{2}\exp\Big(-\frac{y}{2}\Big)\Big)\Big| \le Cn^{-1/2}(\log n)^{5/2}.$$

The proof is now complete by (4.8) and the above two inequalities. □

Proof of Theorem 1.2. We assume that $\mathsf{E}X_{1,1} = 0$, $\mathsf{E}X^2_{1,1} = 1$. The proof of (1.1) is similar to that of (1.8), hence, is omitted. Now we prove (1.10). In view of Lemma 6.1 and the proofs of Theorem 1.1, we have

$$\begin{aligned}\exp\Big(-&\frac{p^2-p}{2}\mathsf{P}(\chi^2(1) \ge l_{n-}(y))\Big) - Cn^{-1/2}(\log n)^{5/2}\\ &\le \mathsf{P}(n\widetilde{L}^2_n - \alpha_p \le y)\\ &\le \exp\Big(-\frac{p^2-p}{2}\mathsf{P}(\chi^2(1) \ge l_{n+}(y))\Big) + Cn^{-1/2}(\log n)^{5/2}.\end{aligned}$$

Moreover, it follows from (4.9), (4.10) and (4.11) that

$$\begin{aligned}\sup_{y\le -2\log_2 n^8}\Big|&\exp\Big(-\frac{p^2-p}{2}\mathsf{P}(\chi^2(1) \ge l_{n\pm}(y))\Big)\\ &- \exp\Big(-\frac{p^2-p}{2}\mathsf{P}(\chi^2(1) \ge \alpha_p + y)\Big)\Big| \le Cn^{-1/2},\end{aligned}$$

$$\begin{aligned}\sup_{y\ge 2\log n}\Big|&\exp\Big(-\frac{p^2-p}{2}\mathsf{P}(\chi^2(1) \ge l_{n\pm}(y))\Big)\\ &- \exp\Big(-\frac{p^2-p}{2}\mathsf{P}(\chi^2(1) \ge \alpha_p + y)\Big)\Big| \le Cn^{-1/2},\end{aligned}$$

and

$$\begin{aligned}\sup_{-2\log_2 n^8 \le y \le 2\log n}\Big|&\exp\Big(-\frac{p^2-p}{2}\mathsf{P}(\chi^2(1) \ge l_{n\pm}(y))\Big)\\ &- \exp\Big(-\frac{p^2-p}{2}\mathsf{P}(\chi^2(1) \ge \alpha_p + y)\Big)\Big| \le Cn^{-1/2}(\log n)^{5/2}.\end{aligned}$$

This completes the proof of Theorem 1.2. □

**5. Proof of the general result.** In this section we outline the proof of the general result, Theorems 3.1 and 3.2.

5.1. *The truncation and notation.* We first truncate $X_{k,i,j}$. Let $\varepsilon$ be a small positive number which will be specified later and put

$$\begin{aligned}
&Y_{k,i,j}^{(m)} = X_{k,i}^{(m)} X_{k,j}^{(m)}, \\
&\widetilde{Y}_{k,i,j}^{(m)} = Y_{k,i,j}^{(m)} I\{|Y_{k,i,j}^{(m)}| \le d^{-1/2}\varepsilon\sqrt{n\log p}\}, \\
&\hat{Y}_{k,i,j}^{(m)} = Y_{k,i,j}^{(m)} I\{|Y_{k,i,j}^{(m)}| \le \sqrt{n}/(\log n)^4\}, \qquad (5.1)\\
&\check{Y}_{k,i,j}^{(m)} = Y_{k,i,j}^{(m)} I\{\sqrt{n}/(\log n)^4 < |Y_{k,i,j}^{(m)}| \le d^{-1/2}\varepsilon\sqrt{n\log p}\}, \\
&\qquad\qquad 1 \le k \le n, 1 \le i, j \le p, 1 \le m \le d
\end{aligned}$$

and

$$\begin{aligned}
&\widetilde{X}_{k,i,j} = (\widetilde{Y}_{k,i,j}^{(1)}, \ldots, \widetilde{Y}_{k,i,j}^{(d)}), \qquad \hat{X}_{k,i,j} = (\hat{Y}_{k,i,j}^{(1)}, \ldots, \hat{Y}_{k,i,j}^{(d)}), \\
&\check{X}_{k,i,j} = (\check{Y}_{k,i,j}^{(1)}, \ldots, \check{Y}_{k,i,j}^{(d)}), \qquad 1 \le k \le n,\ 1 \le i, j \le p, \qquad (5.2)\\
&\widetilde{T}_{p,n} = \max_{1\le i<j\le p}\left\|\sum_{k=1}^{n}\widetilde{X}_{k,i,j}\right\|, \qquad \hat{T}_{p,n} = \max_{1\le i<j\le p}\left\|\sum_{k=1}^{n}\hat{X}_{k,i,j}\right\|, \\
&\widetilde{\mathcal{L}}_{p,n}^2 = \max_{1\le i<j\le p}\frac{\|\sum_{k=1}^{n}\widetilde{X}_{k,i,j}\|^2}{Q_{n,i,j}}.
\end{aligned}$$

5.2. *Auxiliary results.* To prove the general result, we first collect some auxiliary results. As in many previous works on the Poisson approximation, we apply Lemma 5.1 below, which is a special case of Theorem 1 of Arratia, Goldstein and Gordon [1].

LEMMA 5.1 (Arratia, Goldstein and Gordon [1]). *Let $\{\eta_\alpha, \alpha \in I\}$ be random variables on an index set $I$. For each $\alpha \in I$, let $B_\alpha$ be a subset of $I$ with $\alpha \in B_\alpha$. For a given $t \in R$, set $\lambda = \sum_{\alpha\in I}\mathsf{P}(\eta_\alpha > t)$. Then*

$$\left|\mathsf{P}\left(\max_{\alpha\in I}\eta_\alpha \le t\right) - e^{-\lambda}\right| \le (1\wedge\lambda^{-1})(b_1 + b_2 + b_3), \qquad (5.3)$$

*where*

$$\begin{aligned}
&b_1 = \sum_{\alpha\in I}\sum_{\beta\in B_\alpha}\mathsf{P}(\eta_\alpha > t)\mathsf{P}(\eta_\beta > t), \\
&b_2 = \sum_{\alpha\in I}\sum_{\beta\in B_\alpha,\beta\ne\alpha}\mathsf{P}(\eta_\alpha > t, \eta_\beta > t), \\
&b_3 = \sum_{\alpha\in I}\mathsf{E}|\mathsf{P}(\eta_\alpha > t|\sigma(\eta_\beta, \beta \notin B_\alpha)) - \mathsf{P}(\eta_\alpha > t)|
\end{aligned}$$

*and $\sigma(\eta_\beta, \beta \notin B_\alpha)$ is the $\sigma$-algebra generated by $\{\eta_\beta, \beta \notin B_\alpha\}$. In particular, if $\eta_\alpha$ is independent of $\{\eta_\beta, \beta \notin B_\alpha\}$ for each $\alpha$, then $b_3 = 0$.*

The following estimates are also essential for our proof.

PROPOSITION 5.1. *Under the conditions of Theorem 3.1, we have*

$$\begin{aligned} &\sup_{-2\log_2 n^\theta \le y \le 2\log n} \left| \mathsf{P}\left( \left\| \sum_{k=1}^n \widetilde{X}_{k,1,2} \right\| \ge \sqrt{n} y_n \right) - \mathsf{P}(\chi^2(d) \ge \alpha_p + y) \right| \\ &\quad \le C \frac{(\log n)^{5/2}}{p^2 n^{1/2}} \mathsf{E}|XX'|^3 I\{|XX'| \le \sqrt{n}/(\log n)^4\} + Cp^{-3}, \end{aligned} \tag{5.4}$$

*where $y_n = \sqrt{(\alpha_p + y)(1 + O(\sqrt{\log n/n}))}$, $\theta = 8\Gamma(d/2)$, and $\Gamma(\cdot)$ is the gamma function.*

PROPOSITION 5.2. *Under the conditions of Theorem 3.1, we have*

$$\begin{aligned} &\sup_{y \le -2\log_2 n^\theta} \left| \mathsf{P}\left( \frac{\widetilde{T}_{p,n}^2}{n} - \alpha_p \le y \right) - \exp\left( -\frac{p^2 - p}{2} \mathsf{P}(\chi^2(d) \ge \alpha_p + y) \right) \right| \\ &\quad \le Cn^{-2} + Cp^{-1+20\sqrt{\varepsilon}}. \end{aligned}$$

PROPOSITION 5.3. *Under the conditions of Theorem 3.1, we have*

$$\begin{aligned} &\sup_{y \ge 2\log n} \left| \mathsf{P}\left( \frac{\widetilde{T}_{p,n}^2}{n} - \alpha_p \le y \right) - \exp\left( -\frac{p^2 - p}{2} \mathsf{P}(\chi^2(d) \ge \alpha_p + y) \right) \right| \\ &\quad \le C \frac{(\log n)^{5/2}}{n^{1/2}} \mathsf{E}|XX'|^3 I\{|XX'| \le \sqrt{n}/(\log n)^4\} + Cp^{-1}. \end{aligned}$$

PROPOSITION 5.4. *Under the conditions of Theorem 3.1, for any sequence $v_n$ satisfying $v_n \sim 2\sqrt{n \log p}$, we have*

$$\begin{aligned} &\mathsf{P}\left( \left\| \sum_{k=1}^n \widetilde{X}_{k,1,2} \right\| \ge v_n, \left\| \sum_{k=1}^n \widetilde{X}_{k,1,3} \right\| \ge v_n \right) \\ &\quad \le Cp^{-4+20\sqrt{\varepsilon}}. \end{aligned}$$

PROPOSITION 5.5. *Under the conditions of Theorem 3.2, we have*

$$\begin{aligned} &\sup_{y \le -2\log_2 n^\theta} \left| \mathsf{P}(d^2 n \widetilde{\mathcal{L}}_{p,n}^2 - \alpha_p \le y) - \exp\left( -\frac{p^2 - p}{2} \mathsf{P}(\chi^2(d) \ge \alpha_p + y) \right) \right| \\ &\quad \le Cn^{-2} + Cp^{-1+20\sqrt{\varepsilon}} \end{aligned} \tag{5.5}$$

*and*

$$\sup_{y\ge 2\log n}\left|\mathsf{P}(d^2 n\widetilde{\mathcal{L}}^2_{p,n}-\alpha_p\le y)-\exp\left(-\frac{p^2-p}{2}\mathsf{P}(\chi^2(d)\ge\alpha_p+y)\right)\right| \tag{5.6}$$
$$\le C\frac{(\log n)^{5/2}}{n^{1/2}}+\tau_n.$$

Proofs of the propositions above will be given in Section 6.

5.3. *Proof of Theorem 3.1.* Clearly,

$$\mathsf{P}(W^2_{p,n}\ne\widetilde{T}^2_{p,n})\le d\mathsf{P}\left(\max_{1\le k\le n}\max_{1\le i<j\le p}|Y^{(1)}_{k,i,j}|\ge d^{-1/2}\varepsilon\sqrt{n\log p}\right) \tag{5.7}$$
$$\le Cnp^2\mathsf{P}(|XX'|\ge d^{-1/2}\varepsilon\sqrt{n\log p}).$$

To prove Theorem 3.1, it suffices to show that

$$\sup_{y\in R}\left|\mathsf{P}\left(\frac{\widetilde{T}^2_{p,n}}{n}-\alpha_p\le y\right)-\exp\left(-\frac{p^2-p}{2}\mathsf{P}(\chi^2(d)\ge\alpha_p+y)\right)\right| \tag{5.8}$$
$$\le Cp^{-1+20\sqrt{\varepsilon}}+C\frac{(\log n)^{5/2}}{n^{1/2}}\mathsf{E}|XX'|^3 I\{|XX'|\le\sqrt{n}/(\log n)^4\}.$$

Let $I=\{(i,j);1\le i<j\le p\}$ in Lemma 5.1. For $\alpha=(i,j)\in I$, set

$$\eta_\alpha=\left\|\sum_{k=1}^n\widetilde{X}_{k,i,j}\right\|$$

and $B_\alpha=\{(k,l)\in I;\ k=i \text{ or } l=j\}$. Put $t=\sqrt{n\alpha_p+ny}$. Noting that $\eta_\alpha$ is independent of $\{\eta_\beta,\beta\notin B_\alpha\}$ for each $\alpha\in I$, we have

$$|\mathsf{P}(\widetilde{T}_{p,n}\le t)-e^{-\lambda_n}|\le b_{1n}+b_{2n}, \tag{5.9}$$

where

$$\lambda_n=(1/2)(p^2-p)\mathsf{P}\left(\left\|\sum_{k=1}^n\widetilde{X}_{k,1,2}\right\|>t\right),$$

$$b_{1n}\le p^3\mathsf{P}^2\left(\left\|\sum_{k=1}^n\widetilde{X}_{k,1,2}\right\|>t\right),$$

$$b_{2n}\le p^3\mathsf{P}\left(\left\|\sum_{k=1}^n\widetilde{X}_{k,1,2}\right\|>t,\left\|\sum_{k=1}^n\widetilde{X}_{k,1,3}\right\|>t\right).$$

Now (5.8) follows from (5.9) and Propositions 5.1–5.4.

5.4. *Proof of Theorem 3.2.* Assume that $\tau_n \le 1$, otherwise (3.3) is trivial. It suffices to show that

$$\sup_{y\in R}\Big|\mathsf{P}(d^2n\widetilde{\mathcal{L}}^2_{p,n}-\alpha_p\le y)-\exp\Big(-\frac{p^2-p}{2}\mathsf{P}(\chi^2(d)\ge\alpha_p+y)\Big)\Big|$$

$$\le Cp^{-1+20\sqrt{\varepsilon}}+C\frac{(\log n)^{5/2}}{n^{1/2}}+\tau_n.$$

Following the proof of Theorem 3.1, by Lemma 5.1 again, we have

$$|\mathsf{P}(d^2n\widetilde{\mathcal{L}}^2_{p,n}-\alpha_p\le y)-e^{-\widetilde{\lambda}_n}|\le\widetilde{b}_{1n}+\widetilde{b}_{2n}, \tag{5.10}$$

where

$$\widetilde{\lambda}_n=\frac{p^2-p}{2}\mathsf{P}\Bigg(\Bigg\|\sum_{k=1}^n\widetilde{X}_{k,1,2}\Bigg\|>\widetilde{t}_{1,2}\Bigg),$$

$$\widetilde{b}_{1n}\le 2p^3\mathsf{P}^2\Bigg(\Bigg\|\sum_{k=1}^n\widetilde{X}_{k,1,2}\Bigg\|>\widetilde{t}_{1,2}\Bigg),$$

$$\widetilde{b}_{2n}\le p(p^2-p)\mathsf{P}\Bigg(\Bigg\|\sum_{k=1}^n\widetilde{X}_{k,1,2}\Bigg\|>\widetilde{t}_{1,2},\Bigg\|\sum_{k=1}^n\widetilde{X}_{k,1,3}\Bigg\|>\widetilde{t}_{1,3}\Bigg)$$

and

$$\widetilde{t}_{i,j}=\sqrt{\frac{Q_{n,i,j}}{d^2n}(\alpha_p+y)}.$$

Let

$$\hat{A}_{n,i}=\sum_{m=1}^d\sum_{k=1}^n(X^{(m)}_{k,i})^2I\Big\{(X^{(m)}_{k,i})^2\le\sqrt{\frac{n}{\log p}}\Big\},$$

$$\hat{Q}_{n,i,j}=\hat{A}_{n,i}\hat{A}_{n,j},\qquad \hat{t}_{i,j}=\sqrt{\frac{\hat{Q}_{n,i,j}}{d^2n}(\alpha_p+y)}.$$

The main idea of the proof is to replace $Q_{n,i,j}$ by $\hat{Q}_{n,i,j}$ and $\hat{Q}_{n,i,j}$ by some nonrandom constants.

We use Lemma 6.1 and let $\delta=1/4$, $M=2$, $\beta=1$, $a=8D_1(1-M\delta)^{-1}\log p$ and $x=32\times 34d\sqrt{n\log p}$ in (6.1). So $\Sigma_{n,x,a}=0$ and

$$\mathsf{P}\Bigg(\Bigg|\frac{\hat{A}_{n,1}}{dn}-\mathsf{E}X^2I\Big\{X^2\le\sqrt{\frac{n}{\log p}}\Big\}\Bigg|\ge 32\times 34\sqrt{\frac{\log p}{n}}\Bigg)\le Cp^{-4}.$$

Let

$$v_n=\Bigg(n(\alpha_p-2\log_2 n^\theta)\Bigg(1-\mathsf{E}X^2I\Big\{X^2>\sqrt{\frac{n}{\log p}}\Big\}-32\times 34\sqrt{\frac{\log p}{n}}\Bigg)\Bigg)^{1/2}.$$

Then $v_n \sim 2\sqrt{n\log p}$. Note that $A_{n,i,j} \geq \hat{A}_{n.i.j}$ and, hence, $Q_{n,i,j} \geq \hat{Q}_{n,i,j}$. By Proposition 5.4,

$$\begin{aligned}
&\sup_{y\geq -2\log_2 n^8} \mathsf{P}\left(\left\|\sum_{k=1}^n \widetilde{X}_{k,1,2}\right\| > \widetilde{t}_{1,2}, \left\|\sum_{k=1}^n \widetilde{X}_{k,1,3}\right\| > \widetilde{t}_{1,3}\right)\\
&\quad\leq \sup_{y\geq -2\log_2 n^8} \mathsf{P}\left(\left\|\sum_{k=1}^n \widetilde{X}_{k,1,2}\right\| > \hat{t}_{1,2}, \left\|\sum_{k=1}^n \widetilde{X}_{k,1,3}\right\| > \hat{t}_{1,3}\right)\\
&\quad\leq \mathsf{P}\left(\left\|\sum_{k=1}^n \widetilde{X}_{k,1,2}\right\| > v_n, \left\|\sum_{k=1}^n \widetilde{X}_{k,1,3}\right\| > v_n\right) + Cp^{-4}\\
&\quad\leq Cp^{-4+20\sqrt{\varepsilon}}.
\end{aligned} \tag{5.11}$$

Moreover, by Proposition 5.1, for $-2\log_2 n^\theta \leq y \leq 2\log n$,

$$\begin{aligned}
&\mathsf{P}\left(\left\|\sum_{k=1}^n \widetilde{X}_{k,1,2}\right\| > \widetilde{t}_{1,2}\right)\\
&\quad\leq \mathsf{P}\left(\left\|\sum_{k=1}^n \widetilde{X}_{k,1,2}\right\| > \hat{t}_{1,2}\right)\\
&\quad\leq \mathsf{P}\left(\left\|\sum_{k=1}^n \widetilde{X}_{k,1,2}\right\| > \sqrt{n(\alpha_p+y)(1-t_n)}\right) + Cp^{-4}\\
&\quad\leq \mathsf{P}(\chi^2(d) \geq \alpha_p + y)| + C\frac{(\log n)^{5/2}}{p^2 n^{1/2}} + Cp^{-3},
\end{aligned}$$

where $t_n = C\sqrt{\log n/n}$. Note that

$$\begin{aligned}
&\mathsf{P}\left(\left\|\sum_{k=1}^n \widetilde{X}_{k,1,2}\right\| > \widetilde{t}_{1,2}\right)\\
&\quad\geq \mathsf{P}\left(\left\|\sum_{k=1}^n \widetilde{X}_{k,1,2}\right\| > \hat{t}_{1,2}\right)\\
&\qquad - 2d\sum_{i=1}^n \mathsf{P}\left(\left\|\sum_{k=1}^n \widetilde{X}_{k,1,2}\right\| > \hat{t}_{1,2}, (X_{i,1}^{(1)})^2 \geq \sqrt{\frac{n}{\log p}}\right).
\end{aligned}$$

For $y \geq -2\log_2 n^\theta$, letting $\beta = 192\varepsilon$, $q = 4\alpha$, $r = 2 + 2\alpha$, $\delta = 12\varepsilon/5$ and $x = 2(1-2\varepsilon)\sqrt{n\log p}$ in (6.2) yields

$$\mathsf{P}\left(\left\|\sum_{k=1}^n \widetilde{X}_{k,1,2}\right\| > \hat{t}_{1,2}, (X_{i,1}^{(1)})^2 \geq \sqrt{\frac{n}{\log p}}\right)$$

$$\le \mathsf{P}\left(\left\|\sum_{k=1}^{n}\widetilde{X}_{k,1,2}\right\| > \sqrt{n(\alpha_p+y)(1-t_n)}, (X_{i,1}^{(1)})^2 \ge \sqrt{\frac{n}{\log p}}\right) + Cp^{-4}$$

$$\le \mathsf{P}\left(\left\|\sum_{k=1,\neq i}^{n}(\widetilde{X}_{k,1,2}-\mathsf{E}\widetilde{X}_{k,1,2})\right\| > 2(1-2\varepsilon)\sqrt{n\log p}, (X_{i,1}^{(1)})^2 \ge \sqrt{\frac{n}{\log p}}\right)$$

$$+ Cp^{-4}$$

$$\le C[p^{-2(1-10\sqrt{\varepsilon})} + p^{-4}] \times \mathsf{P}\left(|X| \ge \frac{n^{1/4}}{(\log p)^{1/4}}\right) + Cp^{-4}$$

$$\le Cp^{-2(1-10\sqrt{\varepsilon})}\mathsf{P}\left(|X| \ge \frac{n^{1/4}}{(\log p)^{1/4}}\right) + Cp^{-4}.$$

For $-2\log_2 n^\theta \le y \le 2\log n$, it follows from Proposition 5.1 that

$$\mathsf{P}\left(\left\|\sum_{k=1}^{n}\widetilde{X}_{k,1,2}\right\| > \hat{t}_{1,2}\right) \ge \mathsf{P}\left(\left\|\sum_{k=1}^{n}\widetilde{X}_{k,1,2}\right\| > \sqrt{n(\alpha_p+y)(1+t_n)}\right) - Cp^{-4}$$

$$\ge \mathsf{P}(\chi^2(d) \ge \alpha_p + y) - C\frac{(\log n)^{5/2}}{p^2 n^{1/2}} - Cp^{-3}.$$

Combing the above inequalities gives

$$\sup_{-2\log_2 n^\theta \le y \le 2\log n}\left|\mathsf{P}\left(\left\|\sum_{k=1}^{n}\widetilde{X}_{k,1,2}\right\| > \widetilde{t}_{1,2}\right) - \mathsf{P}(\chi^2(d) \ge \alpha_p + y)\right|$$

$$\le C\frac{(\log n)^{5/2}}{p^2 n^{1/2}} + Cp^{-3} \tag{5.12}$$

$$+ Cnp^{-2(1-10\sqrt{\varepsilon})}\mathsf{P}\left(|X| \ge \frac{n^{1/4}}{(\log p)^{1/4}}\right).$$

The proof of Theorem 3.2 is now complete by (5.10)–(5.12) and Proposition 5.5.

**6. Proofs of auxiliary results.** To prove Propositions 5.1–5.4, we need some preliminary lemmas.

The first is a Fuk–Nagaev type inequality for vector-valued random variables.

LEMMA 6.1. *Let* $\xi_i$, $1 \le i \le n$, *be independent random vectors in* $R^d$ *with* $\mathsf{E}\xi_i = 0$ *and* $\mathsf{E}\|\xi_i\|^2 < \infty$, $1 \le i \le n$. *Put* $S_n = \sum_{i=1}^n \xi_i$. *Then for* $0 < \delta < 1$, $\beta > 0$, $a > \delta^{-1}$ *and any* $x > 0$,

$$\mathsf{P}\left(\max_{1\le k\le n}\|S_k\| \ge x + 3\mathsf{E}\|S_n\| + 8\frac{a}{x}\Sigma_{n,x,a}\right)$$

$$\begin{aligned} &\le \sum_{k=1}^{n} \mathsf{P}(\|\xi_k\| > \delta x) + C\left(\frac{\Sigma_{n,x,a}}{x^2}\right)^M \\ &\quad + \exp\left(-\frac{((1-M\delta)x)^2}{2(1+\beta)\Lambda_n}\right) \\ &\quad + \exp\left(-\frac{(1-M\delta)a}{2D_\beta}\right), \end{aligned} \tag{6.1}$$

*where* $\Lambda_n = \sup\{\sum_{k=1}^n \mathsf{E}(u,\xi_k)^2 : \|u\| \le 1\}$, $(\cdot,\cdot)$ *denotes the Euclidean inner product,* $\Sigma_{n,x,a} = \sum_{k=1}^n \mathsf{E}\|\xi_k\|^2 I\{\|\xi_k\| \ge x/a\}$, $D_\beta = 11(1+2/\beta)$, $M$ *is a positive number satisfying* $M\delta < 1$, *and* $C$ *is a constant which depends on* $M$ *and* $\delta$.

*In particular, if* $\max_{1\le k\le n} \mathsf{E}\|\xi_i\|^r \le K$ *for some* $r>2$ *and* $K<\infty$, *then for any* $q \ge 2$ *and* $0 < \beta \le 1$, *there exist* $C_1$, $C_2$ *depending only on* $\beta$, $q$, $K$ *such that, for any* $x \ge C_2\sqrt{n}$ *and* $0 < \delta \le \beta(r-2)(32q+16r-32)^{-1}$,

$$\begin{aligned} &\mathsf{P}\left(\max_{1\le k\le n} \|S_k\| \ge x\right) \\ &\quad \le \sum_{k=1}^{n} \mathsf{P}(\|\xi_k\| > \delta x) + \exp\left(-\frac{x^2}{2(1+\beta)\Lambda_n}\right) + C_1 n^{-q}. \end{aligned} \tag{6.2}$$

PROOF. Put

$$\widetilde{\xi}_i = \xi_i I\{\|\xi_i\| \le \delta x\}, \qquad \widetilde{S}_n = \sum_{i=1}^n \widetilde{\xi}_i,$$

$$\hat{\xi}_i = \xi_i I\{\|\xi_i\| \le x/a\}, \qquad \hat{S}_n = \sum_{i=1}^n \hat{\xi}_i,$$

$$\breve{\xi}_i = \xi_i I\{x/a < \|\xi_i\| \le \delta x\}, \qquad \breve{S}_n = \sum_{i=1}^n \breve{\xi}_i$$

and write

$$B_n = 3\mathsf{E}\|S_n\| + 8\frac{a}{x}\sum_{k=1}^n \mathsf{E}\|\xi_k\|^2 I\{\|\xi_k\| \ge x/a\}.$$

Then

$$\begin{aligned} &\mathsf{P}\left(\max_{1\le k\le n} \|S_k\| \ge x + B_n\right) \\ &\quad \le \mathsf{P}\left(\max_{1\le k\le n} \|\widetilde{S}_k\| \ge x + B_n\right) + \sum_{k=1}^n \mathsf{P}(\|\xi_k\| > \delta x) \end{aligned}$$

(6.3)

$$\le \mathsf{P}\Big(\max_{1\le k\le n}\|\hat S_k\| \ge (1-M\delta)x + B_n/2\Big) + \sum_{k=1}^n \mathsf{P}(\|\xi_k\| > \delta x)$$
$$+ \mathsf{P}\Big(\max_{1\le k\le n}\|\breve S_k\| \ge M\delta x + B_n/2\Big).$$

Since $\mathsf{E}\xi_k = 0$ for $1 \le k \le n$, we have

$$\max_{1\le k\le n}\|\mathsf{E}\hat S_k\| + \tfrac{3}{2}\mathsf{E}\|\hat S_n - \mathsf{E}\hat S_n\| \le B_n/2.$$

Hence,

$$\begin{aligned}
&\mathsf{P}\Big(\max_{1\le k\le n}\|\hat S_k\| \ge (1-M\delta)x + B_n/2\Big)\\
&\quad\le \mathsf{P}\Big(\max_{1\le k\le n}\|\hat S_k - \mathsf{E}\hat S_k\| \ge (1-M\delta)x + \tfrac{3}{2}\mathsf{E}\|\hat S_n - \mathsf{E}\hat S_n\|\Big) \qquad (6.4)\\
&\quad\le \exp\Big(-\frac{((1-M\delta)x)^2}{2(1+\beta)\Lambda_n}\Big) + \exp\Big(-\frac{(1-M\delta)a}{2D_\beta}\Big),
\end{aligned}$$

where the last inequality follows from (3.4) in Einmahl and Li [10].

We now estimate $\mathsf{P}(\max_{1\le k\le n}\|\breve S_k\| \ge M\delta x + B_n/2)$. It follows from the Hoeffding–Bennett inequality that

$$\begin{aligned}
&\mathsf{P}\Big(\max_{1\le k\le n}\|\breve S_k\| \ge M\delta x + B_n/2\Big)\\
&\quad\le \mathsf{P}\Big(\sum_{k=1}^n \|\breve\xi_i\| \ge M\delta x + B_n/2\Big)\\
&\quad\le \mathsf{P}\Big(\sum_{k=1}^n (\|\breve\xi_i\| - \mathsf{E}\|\breve\xi_i\|) \ge M\delta x\Big)\\
&\quad\le \Big(\frac{3\sum_{k=1}^n \mathsf{E}\|\xi_k\|^2 I\{\|\xi_k\| \ge x/a\}}{M\delta^2x^2}\Big)^M.
\end{aligned}$$

So (6.1) is proved.

In order to prove (6.2), we let $a = \max(2D_\tau q(1-M\delta_1)^{-1}\log n, \delta_1^{-1}+1)$, where $\delta_1$ and $\tau$ are positive numbers which will be specified later. Then

$$\exp\Big(-\frac{(1-M\delta_1)a}{2D_\tau}\Big) \le n^{-q}.$$

Note that $\mathsf{E}\|S_n\| \le C\sqrt{n}K^{1/r}$. Since $\Sigma_{n,y,a} \le Kn(ay^{-1})^{r-2}$, we have, for $y \ge \sqrt{n}$, $ay^{-1}\Sigma_{n,y,a} \le Cn^{(3-r)/2}(\log n)^{r-1} \le C_3\sqrt{n}$, and $y^{-2}\Sigma_{n,y,a} \le Cn^{1-r/2} \times a^{r-2} \le Cn^{1-r/2}(\log n)^{r-2}$. Now take $M = (r-2)^{-1}(2q+r-2)$ and $0 < \delta_1 \le$ $\beta(r-2)(16q+8r-16)^{-1}$ such that $M(r/2-1)>q$ and $M\delta_1 \le \beta/8$. Let $C_4 = 8\beta^{-1}(3CK^{1/r}+8C_3)$. By (6.1), it holds that, for $y \ge C_4\sqrt{n}$,

$$\begin{aligned}
&\mathsf{P}\Big(\max_{1\le k\le n}\|S_k\| \ge (1+\beta/8)y\Big)\\
&\quad\le \mathsf{P}\Big(\max_{1\le k\le n}\|S_k\| \ge y + 3\mathsf{E}\|S_n\| + 8\frac{a}{y}\Sigma_{n,y,a}\Big)\\
&\quad\le \sum_{k=1}^{n}\mathsf{P}(\|\xi_k\| > \delta_1 y) + \exp\Big(-\frac{(1-\beta/8)^2y^2}{2(1+\tau)\Lambda_n}\Big) + C_1n^{-q},
\end{aligned}$$

where we let $\tau$ satisfy $(1-\beta/8)^2(1+\tau)^{-1}(1+\beta/8)^{-2} > (1+\beta)^{-1}$. Setting $C_2 = (1+\beta/8)C_4$, $\delta = (1+\beta/8)^{-1}\delta_1$ and $x = y(1+\beta/8)$, we obtain (6.2). $\square$

The following moderate deviation for independent random variables will play an important role in our proof.

Lemma 6.2. *Let $\xi_i$, $1 \le i \le n$, be independent random variables with $\mathsf{E}\xi_i = 0$. Put*

$$s_n^2 = \sum_{i=1}^{n}\mathsf{E}\xi_i^2, \qquad \tau_n = \sum_{i=1}^{n}\mathsf{E}|\xi_i|^3, \qquad S_n = \sum_{i=1}^{n}\xi_i.$$

*Assume that*

$$|\xi_i| \le c_n s_n$$

*for $1 \le i \le n$ and some $0 < c_n \le 1$. Then*

$$\mathsf{P}(S_n \ge xs_n) = e^{\gamma(x/s_n)}(1-\Phi(x))(1+\theta_{n,x}(1+x)s_n^{-3}\tau_n) \tag{6.5}$$

*for $0 < x \le 1/(18c_n)$, where $|\theta_{n,x}| \le 36$ and $\gamma(x)$ is the Cramér–Petrov series (see Petrov [18]) satisfying $|\gamma(x)| \le 2x^3\tau_n s_n^{-3}$. In particular, we have*

$$\mathsf{P}(S_n \ge xs_n) = (1-\Phi(x))(1+\theta_{n,x}(1+x)^3s_n^{-3}\tau_n) \tag{6.6}$$

*for $0 < x \le 1/(18c_n^{1/3})$, where $|\theta_{n,x}| \le 40$.*

Proof. (6.5) is due to Sakhanenko [19] (see Example 1 in Sakhanenko [19]) and (6.6) follows from (6.5). $\square$

The next lemma is simple but useful to provide a moderate deviation of the convolution of independent sums.

LEMMA 6.3. *Let $U_1$, $U_2$, $V_1$ and $V_2$ be independent random variables. Assume that there exist $0 \le c_0 \le 1$ and $x_0$ such that, for any $0 \le x \le x_0$,*

$$\mathsf{P}(|U_1| \ge x) = \mathsf{P}(|V_1| \ge x)(1 + \theta_{1,x}) \tag{6.7}$$

*and*

$$\mathsf{P}(|U_2| \ge x) = \mathsf{P}(|V_2| \ge x)(1 + \theta_{2,x}), \tag{6.8}$$

*where $|\theta_{1,x}| \le c_0$ and $|\theta_{2,x}| \le c_0$. Then*

$$\mathsf{P}(U_1^2 + U_2^2 \ge x^2) = \mathsf{P}(V_1^2 + V_2^2 \ge x^2)(1 + \theta_x) \tag{6.9}$$

*for $0 \le x \le x_0$, where $|\theta_x| \le 3c_0$.*

PROOF. Observe that $\mathsf{P}(|U_1| \ge x) = 1 = \mathsf{P}(|V_1| \ge x)$ for $x \le 0$, so (6.7) and (6.9) remain valid for $x < 0$ with $\theta_{1,x} = 0 = \theta_{2,x}$. Hence, for $0 \le x \le x_0$,

$$\begin{aligned}
\mathsf{P}(U_1^2 + U_2^2 \ge x^2) &= \mathsf{E}\{\mathsf{P}(U_1^2 \ge x^2 - U_2^2 | U_2)\} \\
&\le \mathsf{E}\{\mathsf{P}(V_1^2 \ge x^2 - U_2^2 | U_2)(1 + \theta_{1,x^2 - U_2^2})\} \\
&\le \mathsf{E}\{\mathsf{P}(V_1^2 \ge x^2 - U_2^2 | U_2)(1 + c_0)\} \\
&\le (1 + c_0)\mathsf{P}(U_2^2 \ge x^2 - V_1^2) \\
&= \mathsf{E}\{\mathsf{P}(U_2^2 \ge x^2 - V_1^2 | V_1)\} \\
&\le (1 + c_0)\mathsf{E}\{\mathsf{P}(V_2^2 \ge x^2 - V_1^2)(1 + c_0)\} \\
&= (1 + c_0)^2 \mathsf{P}(V_1^2 + V_2^2 \ge x^2) \\
&\le (1 + 3c_0)\mathsf{P}(V_1^2 + V_2^2 \ge x^2).
\end{aligned} \tag{6.10}$$

Similarly,

$$\begin{aligned}
\mathsf{P}(U_1^2 + U_2^2 \ge x^2) &\ge (1 - c_0)\mathsf{P}(U_2^2 \ge x^2 - V_1^2) \\
&\ge (1 - c_0)^2 \mathsf{P}(V_1^2 + V_2^2 \ge x^2) \\
&\ge (1 - 2c_0)\mathsf{P}(V_1^2 + V_2^2 \ge x^2).
\end{aligned} \tag{6.11}$$

This proves (6.9) by (6.10) and (6.11). □

REMARK 6.1. It is easy to see that Lemma 6.3 remains valid for $m$ independent squared variables.

LEMMA 6.4. *If condition (1.7) is satisfied, then*

$$\mathsf{E}(|X_{1,1}|^{2+4\alpha} / (1 + \log|X_{1,1}|)^{4+4\alpha}) < \infty. \tag{6.12}$$

PROOF. It is easy to see that (1.7) implies

$$\mathsf{E}(|X_{1,1}X_{1,2}|^{2+4\alpha}/(1+\log|X_{1,1}X_{1,2}|)^{4+4\alpha})<\infty,$$

which yields (6.12) by the independence of $X_{1,1}$ and $X_{1,2}$. □

PROOF OF PROPOSITION 5.1. Let $\widetilde{X}_{k,1,2}$, $\hat{X}_{k,1,2}$, $Y^{(1)}_{k,1,2}$ be defined in (5.1) and (5.2). Let

$$A=\bigcup_{i=1}^{d}\bigcup_{k=1}^{n}\left\{|Y^{(i)}_{k,1,2}|>\frac{\sqrt{n}}{(\log n)^4}\right\}.$$

Then

$$\begin{aligned}
&\mathsf{P}\left(\left\|\sum_{k=1}^{n}\widetilde{X}_{k,1,2}\right\|\ge\sqrt{n}y_n\right)\\
&\quad=\mathsf{P}\left(\left\|\sum_{k=1}^{n}\widetilde{X}_{k,1,2}\right\|\ge\sqrt{n}y_n,A\right)+\mathsf{P}\left(\left\|\sum_{k=1}^{n}\widetilde{X}_{k,1,2}\right\|\ge\sqrt{n}y_n,A^c\right)\qquad(6.13)\\
&\quad\le\mathsf{P}\left(\left\|\sum_{k=1}^{n}\widetilde{X}_{k,1,2}\right\|\ge\sqrt{n}y_n,A\right)+\mathsf{P}\left(\left\|\sum_{k=1}^{n}\hat{X}_{k,1,2}\right\|\ge\sqrt{n}y_n\right).
\end{aligned}$$

First, we prove that $\mathsf{P}(\|\sum_{k=1}^{n}\widetilde{X}_{k,1,2}\|\ge\sqrt{n}y_n,A)$ is small, that is,

$$\sup_{-2\log_2 n^\theta\le y\le 2\log n}\mathsf{P}\left(\left\|\sum_{k=1}^{n}\widetilde{X}_{k,1,2}\right\|\ge\sqrt{n}y_n,A\right)\le Cp^{-3}.\qquad(6.14)$$

Noting that for $-2\log_2 n^\theta\le y\le 2\log n$

$$\begin{aligned}
&\mathsf{P}\left(\left\|\sum_{k=1}^{n}\widetilde{X}_{k,1,2}\right\|\ge\sqrt{n}y_n,A\right)\\
&\quad\le d\sum_{i=1}^{n}\mathsf{P}\left(\left\|\sum_{k=1}^{n}\widetilde{X}_{k,1,2}\right\|\ge\sqrt{n}y_n,|Y^{(1)}_{i,1,2}|>\frac{\sqrt{n}}{(\log n)^4}\right)\\
&\quad\le d\sum_{i=1}^{n}\mathsf{P}\left(\left\|\sum_{k=1,k\ne i}^{n}\widetilde{X}_{k,1,2}\right\|\ge 2(1-\varepsilon)\sqrt{n\log p},|Y^{(1)}_{i,1,2}|>\frac{\sqrt{n}}{(\log n)^4}\right)\qquad(6.15)\\
&\quad=d\sum_{i=1}^{n}\mathsf{P}\left(\left\|\sum_{k=1,k\ne i}^{n}\widetilde{X}_{k,1,2}\right\|\ge 2(1-\varepsilon)\sqrt{n\log p}\right)\\
&\qquad\times\mathsf{P}\left(|Y^{(1)}_{i,1,2}|>\frac{\sqrt{n}}{(\log n)^4}\right),
\end{aligned}$$

we have, by Lemma 6.1,

$$
\begin{aligned}
&\mathsf{P}\Bigg(\Bigg\|\sum_{k=1,k\neq i}^{n}\widetilde{X}_{k,1,2}\Bigg\|\geq 2(1-\varepsilon)\sqrt{n\log p}\Bigg)\\
&\quad\leq \mathsf{P}\Bigg(\Bigg\|\sum_{k=1,k\neq i}^{n}X_{k,1,2}\Bigg\|\geq 2(1-\varepsilon)\sqrt{n\log p}\Bigg)\\
&\qquad+dn\mathsf{P}(|Y_{1,1,2}^{(1)}|\geq d^{-1/2}\varepsilon\sqrt{n\log p})\\
&\quad\leq dn\mathsf{P}(|Y_{1,1,2}^{(1)}|\geq d^{-1/2}\varepsilon\sqrt{n\log p})+n\mathsf{P}(\|X_{1,1,2}\|\geq\delta\sqrt{n\log p})\\
&\qquad+\exp\bigg(-\frac{4(1-\varepsilon)^2 n\log p}{2(1+\beta)\Lambda_n}\bigg)+Cn^{-q}\\
&\quad\leq Cn\mathsf{P}(|Y_{1,1,2}^{(1)}|\geq\delta'\sqrt{n\log p})+\exp\bigg(-\frac{2(1-\varepsilon)^2 n\log p}{(1+\beta)\Lambda_n}\bigg)+Cn^{-q},
\end{aligned}
$$

where $\delta$ and $\delta'$ are some positive numbers, $\beta$ is any positive number, $q$ is a large number and $\Lambda_n=(n-1)\sup\{\mathsf{E}(u,X_{1,1,2})^2:\|u\|\leq 1\}=n-1$. For $-2\log_2 n^{\theta}\leq y\leq 2\log n$, by letting $\beta$ small enough such that $(1-\varepsilon)^2(1+\beta)^{-1}>1-2\sqrt{\varepsilon}$, we have

$$
\exp\bigg(-\frac{2(1-\varepsilon)^2 n\log p}{(1+\beta)\Lambda_n}\bigg)\leq Cp^{-2(1-2\varepsilon)^2/(1+\beta)}\leq Cp^{-2(1-2\sqrt{\varepsilon})}.
$$

By Markov's inequality and Lemma 6.4,

$$
np^{-2(1-2\sqrt{\varepsilon})}\mathsf{P}\bigg(|Y_{1,1,2}^{(1)}|>\frac{\sqrt{n}}{(\log n)^4}\bigg)\leq Cp^{-3}
$$

and

$$
n^2\mathsf{P}(|Y_{1,1,2}^{(1)}|\geq\delta'\sqrt{n\log p})\mathsf{P}\bigg(|Y_{1,1,2}^{(1)}|>\frac{\sqrt{n}}{(\log n)^4}\bigg)\leq Cp^{-3}.
$$

This proves (6.14).

Now we come to estimate $\mathsf{P}(\|\sum_{k=1}^n\hat{X}_{k,1,2}\|\geq\sqrt{n}y_n)$. Observe that

$$
\begin{aligned}
&\mathsf{P}\Bigg(\Bigg\|\sum_{k=1}^{n}\hat{X}_{k,1,2}\Bigg\|\geq\sqrt{n}y_n\Bigg)\\
&\quad\leq \mathsf{P}\Bigg(\Bigg\|b_n^{-1/2}\sum_{k=1}^{n}(\hat{X}_{k,1,2}-\mathsf{E}\hat{X}_{k,1,2})\Bigg\|\geq b_n^{-1/2}(\sqrt{n}y_n-dna_n)\Bigg) \qquad (6.16)\\
&\quad\leq \mathsf{P}\Bigg(\Bigg\|(nb_n)^{-1/2}\sum_{k=1}^{n}(\hat{X}_{k,1,2}-\mathsf{E}\hat{X}_{k,1,2})\Bigg\|\geq y_n-c_n\Bigg),
\end{aligned}
$$

where $a_n = \mathsf{E}\|X_{1,1,2}\| I\{\|X_{1,1,2}\| \geq \sqrt{n}/(\log n)^4\}$, $b_n = \mathsf{Var}(\hat{Y}^{(1)}_{1,1,2})$ and

$$\begin{aligned} c_n &= C\sqrt{n} a_n + C\sqrt{\log n}\mathsf{E}|Y^{(1)}_{1,1,2}|^2 I\{|Y^{(1)}_{1,1,2}| \geq \sqrt{n}/(\log n)^4\} \\ &\leq C_\delta p^{-2+\delta} \end{aligned}$$

for any $\delta > 0$. By Lemma 6.2, for $0 \leq x \leq (\log n)^{4/3}/100$, we have

$$\begin{aligned} &\mathsf{P}\left(\left|b_n^{-1/2}\sum_{k=1}^n (\hat{Y}^{(1)}_{k,1,2} - \mathsf{E}\hat{Y}^{(1)}_{k,1,2})\right| \geq \sqrt{n}x\right) \\ &\qquad = 2(1-\Phi(x))(1+\theta_{n,x}(1+x)^3(nb_n)^{-3/2}\tau_n), \end{aligned}$$

where $\theta_{n,x} \leq 40 \times 8$ and $\tau_n = n\mathsf{E}|\hat{Y}^{(1)}_{1,1,2}|^3$. So by Lemma 6.3 (see also Remark 6.1), we get, for $0 \leq x \leq (\log n)^{4/3}/100$,

$$\begin{aligned} &\mathsf{P}\left(\left\|b_n^{-1/2}\sum_{k=1}^n (\hat{X}_{k,1,2} - \mathsf{E}\hat{X}_{k,1,2})\right\|^2 \geq nx^2\right) \\ &\qquad = \mathsf{P}(\chi^2(d) \geq x^2)(1+\theta'_{n,x}(1+x)^3(nb_n)^{-3/2}\tau_n), \end{aligned} \tag{6.17}$$

where $\theta'_{n,x} \leq 3^{d-1}320$. Putting (6.16) with (6.17) together yields, for $-2\log_2 n^\theta \leq y \leq 2\log n$ and some $0 < C < \infty$ (not depending on $y$),

$$\begin{aligned} &\mathsf{P}\left(\left\|\sum_{k=1}^n \hat{X}_{k,1,2}\right\| \geq \sqrt{n}y_n\right) \\ &\quad \leq \mathsf{P}(\chi^2(d) \geq y_n^2) + c_n y_n^{d-1}\exp(-(y_n - c_n)^2/2) \\ &\qquad + C\mathsf{P}(\chi^2(d) \geq (y_n - c_n)^2)(1+y_n)^3(nb_n)^{-3/2}\tau_n \\ &\quad \leq \mathsf{P}(\chi^2(d) \geq \alpha_p + y) + Cp^{-2}n^{-1/2}(\log n)^{5/2} \\ &\qquad + C(p^{-2}c_n(\log n)^{3/2} + p^{-2}(\log n)^{5/2}n^{-1/2}\mathsf{E}|\hat{Y}^{(1)}_{1,1,2}|^3) \\ &\quad \leq \mathsf{P}(\chi^2(d) \geq \alpha_p + y) + C(p^{-3} + p^{-2}(\log n)^{5/2}n^{-1/2}\mathsf{E}|\hat{Y}^{(1)}_{1,1,2}|^3). \end{aligned}$$

Similarly, we have, for $-2\log_2 n^\theta \leq y \leq 2\log n$,

$$\begin{aligned} &\mathsf{P}\left(\left\|\sum_{k=1}^n \hat{X}_{k,1,2}\right\| \geq \sqrt{n}y_n\right) \\ &\quad \geq \mathsf{P}(\chi^2(d) \geq \alpha_p + y) - C(p^{-3} + p^{-2}(\log n)^{5/2}n^{-1/2}\mathsf{E}|\hat{Y}^{(1)}_{1,1,2}|^3). \end{aligned}$$

This completes the proof of Proposition 5.1 by combining the above inequalities. $\square$

Proof of Proposition 5.2. Since

$$\mathsf{P}(\chi^2(d) \geq x) \sim 2^{1-d/2}\Gamma^{-1}(d/2)x^{d/2-1}\exp(-x/2) \tag{6.18}$$

as $x \to \infty$, we have

$$\sup_{y \leq -2\log_2 n^\theta} \exp\biggl(-\frac{p^2-p}{2}\mathsf{P}(\chi^2(d) \geq \alpha_p + y)\biggr) \leq Cn^{-2}. \tag{6.19}$$

Noting that

$$\sup_{y \leq -2\log_2 n^\theta} \mathsf{P}(\widetilde{T}_{p,n}^2 - n\alpha_p \leq ny) \leq \mathsf{P}(\widetilde{T}_{p,n}^2 \leq n\alpha_p - 2n\log_2 n^\theta), \tag{6.20}$$

Propositions 5.1 and 5.4 imply

$$\mathsf{P}(\widetilde{T}_{p,n}^2 \leq n\alpha_p - 2n\log_2 n^\theta) \leq e^{-\varphi_n} + Cp^{-1+20\sqrt{\varepsilon}}, \tag{6.21}$$

where

$$\varphi_n = \frac{p^2-p}{2}\mathsf{P}\Biggl(\Biggl\|\sum_{k=1}^n \widetilde{X}_{k,1,2}\Biggr\| \geq \sqrt{n\alpha_p - 2n\log_2 n^\theta}\Biggr).$$

By Proposition 5.1 again, we have $\varphi_n \geq 2\log n$ for $n$ large. Now Proposition 5.2 follows by (6.19), (6.20) and (6.21). □

Proof of Proposition 5.3. By (6.18),

$$\sup_{y \geq 2\log n}\biggl|1 - \exp\biggl(-\frac{p^2-p}{2}\mathsf{P}(\chi^2(d) \geq \alpha_p + y)\biggr)\biggr| \leq Cn^{-1}. \tag{6.22}$$

Moreover, by Proposition 5.1,

$$\begin{aligned}
&\sup_{y \geq 2\log n} \mathsf{P}(\widetilde{T}_{p,n}^2 \\
&\qquad \geq n\alpha_p + ny) \leq \mathsf{P}(\widetilde{T}_{p,n}^2 \geq n\alpha_p + 2n\log n) \\
&\qquad \leq p^2\mathsf{P}\Biggl(\Biggl\|\sum_{k=1}^n \widetilde{X}_{k,1,2}\Biggr\|^2 \geq n\alpha_p + 2n\log n\Biggr) \\
&\qquad \leq C\frac{(\log n)^{5/2}}{n^{1/2}}\mathsf{E}|XX'|^3 I\{|XX'| \leq \sqrt{n}/(\log n)^4\} \\
&\qquad\quad + Cp^{-1} + Cn^{-1}.
\end{aligned} \tag{6.23}$$

Proposition 5.3 now follows from (6.22) and (6.23). □

PROOF OF PROPOSITION 5.4. Observe that

$$
\begin{aligned}
&\mathsf{P}\left(\left\|\sum_{i=1}^{n}\widetilde{X}_{i,1,2}\right\| \ge v_n, \left\|\sum_{i=1}^{n}\widetilde{X}_{i,1,3}\right\| \ge v_n\right)\\
&\quad \le \mathsf{P}\left(\left\|\sum_{i=1}^{n}(\widetilde{X}_{i,1,2} - \mathsf{E}\widetilde{X}_{i,1,2})\right\| \ge (1-\varepsilon)v_n,\right.\\
&\qquad\qquad \left.\left\|\sum_{i=1}^{n}(\widetilde{X}_{i,1,3} - \mathsf{E}\widetilde{X}_{i,1,3})\right\| \ge (1-\varepsilon)v_n\right)\\
&\quad \le \mathsf{P}\left(\left\|\sum_{i=1}^{n}(\widetilde{X}_{i,1,2} - \mathsf{E}\widetilde{X}_{i,1,2}, \widetilde{X}_{i,1,3} - \mathsf{E}\widetilde{X}_{i,1,3})\right\|^2 \ge 2((1-\varepsilon)v_n)^2\right)\\
&\quad \le \mathsf{P}\left(\left\|\sum_{i=1}^{n}(\widetilde{X}_{i,1,2} - \mathsf{E}\widetilde{X}_{i,1,2}, \widetilde{X}_{i,1,3} - \mathsf{E}\widetilde{X}_{i,1,3})\right\|^2\right.\\
&\qquad\qquad\qquad \left.\ge 2((2-3\varepsilon)\sqrt{n\log p})^2\right).
\end{aligned} \tag{6.24}
$$

Take $\beta = 192\varepsilon$, $q = 4\alpha$, $r = 2 + 2\alpha$, $\delta = \beta(r-2)(32q + 16r - 32)^{-1} = 12\varepsilon/5$ and $x = \sqrt{2}(2-3\varepsilon)\sqrt{n\log p}$ in (6.2). Since $\|\widetilde{X}_{i,1,2}\| \le \varepsilon\sqrt{n\log p}$, we have

$$
\mathsf{P}(\|(\widetilde{X}_{i,1,2} - \mathsf{E}\widetilde{X}_{i,1,2}, \widetilde{X}_{i,1,3} - \mathsf{E}\widetilde{X}_{i,1,3})\| \ge \delta x) = 0. \tag{6.25}
$$

It can be shown that the largest eigenvalue of the covariance matrix of $(\widetilde{X}_{1,1,2} - \mathsf{E}\widetilde{X}_{1,1,2}, \widetilde{X}_{1,1,3} - \mathsf{E}\widetilde{X}_{1,1,3})$ tends to 1 as $n \to \infty$. Therefore, for $0 < \varepsilon < 10^{-4}$,

$$
\exp\left(-\frac{x^2}{2(1+192\varepsilon)\Lambda_n}\right) \le C\exp\left(-\frac{(2-3\varepsilon)^2\log p}{1+193\varepsilon}\right) \le Cp^{-4+20\sqrt{\varepsilon}}. \tag{6.26}
$$

Proposition 5.4 is proved by (6.2) and (6.24)–(6.26). □

PROOF OF PROPOSITION 5.5. The proof of (5.5) follows from the proof of Proposition 5.2, while the proof of (5.6) is similar to that of Proposition 5.3. The details are omitted here. □

PROOF OF (1.11). Let

$$
\lambda_n = \frac{p^2 - p}{2}\mathsf{P}(\chi^2(1) \ge y_n^2), \qquad y_n = \sqrt{4\log p - \log_2 p + y}.
$$

It is known that

$$
\frac{2}{\sqrt{2\pi}y_n}\left(1 - \frac{1}{y_n^2}\right)e^{-y_n^2/2} \le \mathsf{P}(\chi^2(1) \ge y_n^2) \le \frac{2}{\sqrt{2\pi}y_n}e^{-y_n^2/2}. \tag{6.27}
$$

Thus,

$$\lambda_n \le \frac{p^2 - p}{\sqrt{2\pi} y_n} \exp\left(-\frac{y_n^2}{2}\right) = \frac{(1 - p^{-1})\sqrt{\log p}}{\sqrt{2\pi} y_n} \exp\left(-\frac{y}{2}\right).$$

By (6.27) again,

$$\lambda_n \ge \frac{(1 - p^{-1})\sqrt{\log p}}{\sqrt{2\pi} y_n} \exp\left(-\frac{y}{2}\right) - O\left(\frac{1}{\log n}\right).$$

Hence,

$$\lambda_n = \frac{1}{\sqrt{8\pi}} \exp\left(-\frac{y}{2}\right) + A_n \frac{1}{\sqrt{8\pi}} \exp\left(-\frac{y}{2}\right) + O\left(\frac{1}{\log n}\right),$$

where $A_n \sim \log_2 n/(8 \log n)$ as $n \to \infty$. Finally, we can write

$$e^{-\lambda_n} = (1 - B_n) \exp\left(-\frac{1}{\sqrt{8\pi}} \exp\left(-\frac{y}{2}\right)\right), \qquad B_n \sim \frac{\log_2 n}{8 \log n} \frac{1}{\sqrt{8\pi}} \exp\left(-\frac{y}{2}\right).$$

This proves (1.11). $\square$

**Acknowledgments.** We thank an anonymous referee and an Associate Editor for their insightful comments.

W.-D. Liu
Department of Mathematics
Zhejiang University
Hangzhou 310027
China
E-mail: liuweidong99@gmail.com

Z. Lin
Department of Mathematics
Zhejiang University
Hangzhou 310027
China
E-mail: zlin@zju.edu.cn

Q.-M. Shao
Department of Mathematics
Hong Kong University of Science
and Technology
Clear Water Bay
Kowloon
Hong Kong
China
E-mail: maqmshao@ust.hk